\documentclass[11pt,reqno]{amsart}
\usepackage{amsmath} 
\usepackage{amssymb}
\usepackage{dsfont}

\def\squarebox#1{\hbox to #1{\hfill\vbox to #1{\vfill}}}

\theoremstyle{plain}

\newtheorem{Thm}{Theorem}

\newtheorem{lem}{Lemma}

\newcommand{\R}{\mathbb{R}}
 
\newcommand{\B}{\dot{\mathcal{H}}_1({\R}^n)}

\newcommand{\CI}{{\mathcal C}^{\infty}_{0}({\R}^{n}) }

\newcommand{\half}{\frac{1}{2}}

\def\phi {\varphi}

\newtheorem{rem}{Remark}

\def\dem{\noindent    {\bf Proof.} }

\numberwithin{equation}{section}

\begin{document}

\def\dem{\noindent    {\bf Proof:} }

\title[Semilinear wave equation ]
{ Cauchy problem for semilinear wave equation with time-dependent metrics}
\author[Y. Kian]{Yavar Kian}

\address {Universit\'e Bordeaux I, Institut de Math\'ematiques de Bordeaux,  351, Cours de la Lib\'eration, 33405  Talence, France}
\email{Yavar.Kian@math.u-bordeaux1.fr}
\maketitle
\begin{abstract}
We establish  the existence of   weak solutions $u$ of the semilinear wave equation\\
 $\partial_t^2 u-\textrm{div}_x(a(t,x)\nabla_xu)=f_k(u)$ where $a(t,x)$ is equal to $1$ outside a compact set with respect to $x$ and a non-linear term $f_k$ which satisfies $\vert f_k(u)\vert\leq C\vert u\vert^k$. For some  non-trapping time-periodic perturbations $a(t,x)$, we obtain the long time existence of  solution for small initial data.
\end{abstract}

\section{Introduction}

Consider the semilinear Cauchy problem
\begin{equation} \label{eq:1}  \left\{\begin{array}{c}
u_{tt}-\textrm{div}_{x}(a(t,x)\nabla_{x}u)-f_k(u)=0,\ \ (t,x)\in{\R}^{n+1},\\
(u,u_t)(0,x)=(g_1(x),g_2(x))=g(x),\ \ x\in{\R}^n,\end{array}\right.\end{equation}
where for a given $k>1$ the  non-linearity $f_k$ is assumed to be a $C^1$ function on ${\R}$ satisfying $f_k(0)=0$, $\vert f'_k(u)\vert\leq C\vert u\vert^{k-1}$ and the perturbation $a(t,x)\in C^\infty({\R}^{n+1})$ satisfies the conditions:
\begin{equation}\label{eq=1}\begin{array}{l}(i) \ C_0\geq a(t,x)\geq c_0>0,\   \forall(t,x)\in{\R}^{n+1},\\
(ii)\ \textrm{ there exists }\rho>0\textrm{ such that }a(t,x)=1\textrm{ for }\vert x\vert\geq\rho.\end{array}\end{equation}
 Denote by $\dot{H}^1({\R}^n)$ the closure of ${\CI}$ with respect to the norm
\[\Vert \phi\Vert_{\dot{H}^1}=\left(\int_{{\R}^n}\vert\nabla\phi(x)\vert^2dx\right)^\half.\]
Throughout this paper we assume that $n\geq3$ and that the initial data $g$ is  in the energy space ${\B}=\dot{H}^1({\R}^n)\times L^2({\R}^n)$. Consider the linear problem associated to (1.1)
\begin{equation} \label{eq:2}  \left\{\begin{array}{c}
u_{tt}-\textrm{div}_{x}(a(t,x)\nabla_{x}u)=0,\ \ (t,x)\in{\R}^{n+1},\\
(u,u_{t})(s,x)=(g_{1}(x),g_{2}(x))=g(x),\ \ x\in{\R}^n,\end{array}\right.\end{equation}
where $g\in{\B}$. The solution of (1.3) is given by the propagator
\[\mathcal{U}(t,s):{\B}\ni(g_1,g_2)=g\mapsto \mathcal{U}(t,s)g=(u,u_t)(t,x)\in{\B}.\]

We denote by $U(t,s)$ and $V(t,s)$ the operators defined by
\[ U(t,s)f=\left(\mathcal{U}(t,s)(f,0)\right)_1,\quad f\in\dot{H}^1({\R}^n),\]
\[ V(t,s)h=\left(\mathcal{U}(t,s)(0,h)\right)_1,\quad h\in L^2({\R}^n),\]
where $(h_1,h_2)_1=h_1$.
We say that $u\in\mathcal{C}([0,T_1],\dot{H}^1)$ is a weak solution of (1.1) if  for all $t\in[0,T_1]$ we have
\begin{equation}\label{eq:1}
\begin{array}{lll}u(t)&=&\left(\mathcal{U}(t,0)g+\int_0^t\mathcal{U}(t,s)(0,f_k(u(s))) \textrm{d}s\right)_1\\

\ &=&\left(\mathcal{U}(t,0)g
\right)_1+\int_0^tV(t,s)( f_k(u(s)) \textrm{d}s.\\ \end{array}
 \end{equation}

Let $a(t,x)=1$. Then we have the following Cauchy problem
\begin{equation} \label{eq:2}  \left\{\begin{array}{c}
u_{tt}-\Delta_x u-f_k(u)=0,\ \ (t,x)\in{\R}^{n+1},\\
(u,u_{t})(0,x)=(g_{1}(x),g_{2}(x))=g(x),\ \ x\in{\R}^n.\end{array}\right.\end{equation}
The problem (1.5) has been extensively studied for $g\in{\B}$. For example, the global well-posedness of the problem (1.5) has been established for the case of the sub-critical growth $1<k<1+\frac{4}{n-2}$ (see [6] and [19]) or for the case of the critical growth $k=1+\frac{4}{n-2}$ (see [15] and [19]). For the case $k>1+\frac{4}{n-2}$ it is not yet clear whether there exists or not a global regular solution for the Cauchy problem (1.5) with arbitrary initial data. On the other hand, local well posedness as well as global well-posedness, with small initial data in fractional Sobolev spaces have been also studied by many authors for the problem (1.5) under minimal regularity assumptions on the initial data (see [8] and [19]).

In [18] Michael Reissig and Karen Yagdjian established Strichartz decay estimates for the solution of strictly hyperbolic equations of second order with coefficients  depending only on $t$. We can apply these estimates to prove existence results for the  solution of problem (1.1) when $a(t,x)=a(t)$ is independent on $x$ (see [11] and [20] for the case of the free wave equation). It seems that our paper is one of the first works where one treats non-linear wave equations with time dependent perturbations $a(t,x)$ depending on $t$ and $x$.

The goal of this paper is to find sufficient conditions for the existence of a weak solution of  (1.1)  when $0\leq t\leq T_1$.
For this purpose, we will use  Strichartz estimates to study local and long time existence and uniqueness of solutions of the problem (1.1).  In fact, for suitable $k$ Strichartz estimates allow us to find a fixed point of the map
\[\mathcal{G}(u)=\left(\mathcal{U}(t,0)g\right)_1+\int_0^tV(t,s)f_k(u(s)) \textrm{d}s,\]
 in $\mathcal{C}([0,T_1],\dot{H}^1)$ for some well chosen $k>1$. The fixed point of $\mathcal{G}$ is  local weak solution of (1.1).
In [9] we have  established  local homogeneous  Strichartz estimates  for $n\geq3$ and $a(t,x)$ satisfying (1.2), and global homogeneous Strichartz estimates  when  $n\geq3$ is odd for some non-trapping time-periodic perturbation  $a(t,x)$ ( see Section 2). Recently global Strichartz estimates for even dimensions $n\geq4$ have been obtained in [10]. One way to obtain global weak solutions is to apply global non homogeneous Strichartz estimates concerning the solution of the Cauchy problem for $u_{tt} - {\rm div}\:(a(t,x) \nabla u(x)) = G(t,x).$ This leads to some difficulties and this case is not covered by our results in [9] and [10]. On the other hand, for time dependent perturbations we have no conservation laws. For these reasons we obtain only long time existence of weak solution in Section 4.
 In Section 2 we recall the estimates for the linear wave equation with metric $a(t,x)$. In Section 3 we obtain local existence results, while in Section 4 we deal with long time existence.

\begin{rem}
Let the metric $(a_{ij}(t,x))_{1\leq i,j\leq n}$ be such that for all $i,j=1\cdots n$ we have
\[\begin{array}{l}
\displaystyle(i)\ \textrm{ there exists }\rho>0\textrm{ such that }a_{ij}(t,x)=\delta_{ij},\textrm{ for }\vert x\vert\geq\rho,\textrm{ with $\delta_{ij}=0$ for $i\neq j$ and $\delta_{ii}=1$},\\
\displaystyle(ii)\textrm{ there exists }T>0\textrm{ such that }a_{ij}(t+T,x)=a_{ij}(t,x),\    (t,x)\in{\R}^{n+1},\\
\displaystyle(iii)a_{ij}(t,x)=a_{ji}(t,x),  (t,x)\in{\R}^{n+1},\\
\displaystyle(iv)\textrm{ there exist } C_0>c_0>0 \textrm{ such that }C_0\vert\xi\vert^2\geq\sum_{i,j=1}^na_{ij}(t,x)\xi_i\xi_j\geq c_0\vert\xi\vert^2 ,\ \  (t,x)\in{\R}^{1+n},\ \xi\in{\R}^n.\\
\end{array}\]

If we replace $a(t,x)$ in \emph{(1.1)} we get the following problem
\[ \left\{\begin{array}{c}
\displaystyle u_{tt}-\sum_{i,j=1}^n\frac{\partial}{\partial x_i}\left(a_{ij}(t,x)\frac{\partial}{\partial x_j}u\right)-f_k(u)=0,\ \ (t,x)\in{\R}^{n+1},\\
\ \\
(u,u_{t})(s,x)=(f_{1}(x),f_{2}(x))=f(x),\ \ x\in{\R}^n.\end{array}\right.\]

Since, with the same conditions as \emph{(1.1)}, global Strichartz  estimates are true for solutions of this equation when $f_k(u)=0$ \emph{(}see \cite{K3}\emph{)}, all the results of this paper remain true for this problem.
\end{rem}

\section{Strichartz estimates for the linear equation}

In this section we recall some results concerning Strichartz estimates for the problem (1.3). We suppose that $a(t,x)$ satisfies
the conditions (1.2). It was established in [9] that  we have the following estimates.

\begin{Thm}
 Assume $n\geq3$ and let $a(t,x)$ be a $C^\infty$ function on ${\R}^{n+1}$ satisfying conditions  (1.2). Let $2\leq p,q<+\infty$, $\gamma>0$ be such that
\begin{equation}\label{eq:1}\frac{1}{p}=\frac{n(q-2)}{2q}-\gamma,\quad\frac{1}{p}\leq\frac{(n-1)(q-2)}{4q}.\end{equation}
Then there exists $\delta>0$ such that for the solution $u(t,x)$ of (1.3) with $s=0$ we have
\begin{equation}\label{eq:1}\Vert u\Vert_{L^p([0,\delta],L^q({\R}_x^n))}+ \Vert u(t)\Vert_{\mathcal{C}([0,\delta],\dot{H}^\gamma({\R}^n_x))}+\Vert \partial_t(u)(t)\Vert_{\mathcal{C}([0,\delta],\dot{H}^{\gamma-1}({\R}^n_x))}\leq C(p,q,\rho,n)\Vert g\Vert_{\dot{\mathcal{H}}_\gamma}.\end{equation}
\end{Thm}

Now,  let $a(t,x)$ be $T$-periodic with respect to $t$ which means
\[a(t+T,x)=a(t,x),\quad\forall(t,x)\in{\R}^{n+1}.\]
 Moreover, we impose two hypothesis.
The first one says that the perturbation $a(t,x)$ is non-trapping. More precisely, consider the null bicharacteristics $(t(\sigma),x(\sigma),\tau(\sigma),\xi(\sigma))$ of the principal symbol $\tau^2-a(t,x)\vert\xi\vert^2$ of $\partial^2_t-\textrm{div}(a\nabla_xu)$ satisfying
\[t(0)=0,\vert x(0)\vert\leq\rho,\quad \tau^2(\sigma)=a(t(\sigma),x(\sigma))\vert\xi(\sigma)\vert^2.\]
We introduce the following condition.

$\ \ $(H1) We say that the metric $a(t,x)$ is non-trapping if for each $R>\rho$ there exists $S_R>0$ such that

$\ \ \ \ \ \ \ $ $\vert x(\sigma)\vert>R$ for $\vert \sigma\vert\geq S_R$.\\
 Notice that if we have trapping metrics, there exist solutions of (1.3) whose local energy is exponentially growing (see [2]). Thus for trapping metrics it is not possible to establish global Strichartz estimates.

Let $\psi_1,\psi_2\in{\CI}$. We define the cut-off resolvent associated to problem (1.3) by $R_{\psi_1,\psi_2}(\theta)=\psi_1(\mathcal{U}(T,0)-e^{-i\theta})^{-1}\psi_2$. Consider the following assumption.

$\ \ $(H2) Let $\psi_1,\psi_2\in{\CI}$ be such that $\psi_i = 1$ for $|x| \leq \rho + 1 + 3T, i = 1, 2$. Then the operator $R_{\psi_1,\psi_2}(\theta)$ admits a holomorphic extension from $\{\theta\in\mathbb{C}\  :\  \textrm{Im}(\theta) \geq A > 0\}$ to
$\{\theta\in\mathbb{C}\  :\  \textrm{Im}(\theta) \geq  0\}$, for $n \geq 3$, odd, and to $\{\theta\in\mathbb{C}\  :\  \textrm{Im}(\theta) >  0\}$ for $n \geq4$, even.
Moreover, for n even, $R_{\psi_1,\psi_2}(\theta)$ admits a continuous extension from $\{\theta\in\mathbb{C}\  :\  \textrm{Im}(\theta) > 0\}$ to $\{\theta\in\mathbb{C}\  :\  \textrm{Im}(\theta) \geq  0,\theta\neq 2k\pi,\forall k\in\mathbb{Z}\}$ and we have
\[\limsup_{\textrm{Im}(\lambda)>0,\lambda\rightarrow0}\Vert R_{\psi_1,\psi_2}(\lambda)\Vert<\infty.\]

Assuming conditions (H1) and (H2), we obtained the following estimates (see [9], [10]).
\begin{Thm}
Assume $n\geq3$ and let $a(t,x)$ be a $T$-periodic metric satisfying (1.2) for which the conditions (H1) and (H2) are fulfilled. Let $2\leq p,q<+\infty$ be such that
\begin{equation}\label{eq:1}p>2,\quad\frac{1}{p}=\frac{n(q-2)}{2q}-1,\quad\frac{1}{p}\leq\frac{(n-1)(q-2)}{4q}.\end{equation}
Then for the solution $u(t)$ of (1.3) with $s=0$ we have for all $t>0$ the  estimate
\begin{equation} \label{eq:2}\Vert u(t)\Vert_{L^p({\R}_t^+,L^q({R}_x^n))}+ \Vert u(t)\Vert_{\dot{H}^1({\R}^n_x)}+\Vert \partial_t(u)(t)\Vert_{L^2({\R}^n_x)}\leq C(p,q,\rho,T)(\Vert g_1\Vert_{\dot{H}^1({\R}^n)}+\Vert g_2\Vert_{L^2({\R}^n)}).\end{equation}
\end{Thm}
The crucial point in the proof of the global estimates (2.4) is the $L^2$ integrability with respect to $t$ of the local energy (see [9], [10]). For this purpose we need to show that for cut-off functions $\psi_1,\psi_2\in{\CI}$ such that $\psi_i = 1$ for $|x| \leq \rho + 1 + 3T, i = 1, 2$, for $t\geq s$ we have
\begin{equation} \label{eq:2}\Vert\psi_1 \mathcal U(t,s)\psi_2\Vert_{\mathcal L({\B})}\leq C_{\psi_1,\psi_2}d(t-s)\end{equation}
with $d(t)\in L^1({\R}^+)$. To obtain (2.5), we use the assumption (H2). For $n\geq3$, odd, we have an exponential decay of energy and $d(t)=e^{-\delta t}$, $\delta>0$. For $n\geq4$, even, we have another decay. In particular, the estimate
\begin{equation} \label{eq:2}\Vert\psi_1 \mathcal U(NT,0)\psi_2\Vert_{\mathcal L({\B})}\leq \frac{C_{\psi_1,\psi_2}}{(N+1)\ln^2(N+e)},\quad\forall N\in\mathbb N,\end{equation}
implies (2.5). On the other hand, if (2.6) holds, the assumption (H2) for $n$ even is fulfilled. Indeed, for large $A>>1$ and $\textrm{Im}(\theta)\geq AT$ we have
\[R_{\psi_1,\psi_2}(\theta)=-e^{i\theta}\sum_{N=0}^{\infty}\psi_1\mathcal U(NT,0)\psi_2e^{iN\theta}\]
and applying (2.6), we conclude that $R_{\psi_1,\psi_2}(\theta)$ admits a holomorphic extension from $\{\theta\in\mathbb{C}\  :\  \textrm{Im}(\theta) \geq A > 0\}$   to $\{\theta\in\mathbb{C}\  :\  \textrm{Im}(\theta) >  0\}$. Moreover, $R_{\psi_1,\psi_2}(\theta)$ is bounded for $\theta\in{\R}$. We refer to [10] for examples of metrics $a(t,x)$ such that  (2.6) is fulfilled. We
like to mention that in the study of the time-periodic perturbations of the Schr\"odinger operators
(see [3] ) the resolvent of the monodromy operator $(\mathcal U(T) - z)^{-1}$ plays a central role. Moreover, the absence of eigenvalues $z \in\mathbb C, |z| = 1$ of $\mathcal U(T)$, and the behavior of the resolvent for $z$ near $1$, are closely related to the decay of local energy as $t \rightarrow\infty$. So our results may be considered as a natural extension of those for Schr\"odinger operator. On the other hand, for the
wave equation we may have poles $\theta \in\mathbb C$, $ Im\theta > 0$ of the $R_{\psi_1,\psi_2}(\theta)$, while for the Schr\"odinger operator
with time-periodic potentials  such a phenomenon is excluded.

\section{Local time existence}
In this section we assume  $n\geq3$ and let  $a(t,x)$ be a $C^\infty$ function on ${\R}^{n+1}$ satisfying the conditions  (1.2). Motivated by  the work of T. Tao and M. Keel in [12], we will apply Theorem 1 to find $k>1$ for which the problem (1.1) is locally well-posed. For this purpose we need to find $k>1$ so that there exist $2\leq p,q<+\infty$ satisfying (2.1) with $\gamma=1$ for which we have
\begin{equation} \label{eq:2} k=\frac{q}{2},\quad\frac{k}{p}<1.\end{equation}
Then it is easy to see that $k>1$ satisfies (3.1) with $p,q$ satisfying (2.1), if  the following conditions are fulfilled:
\begin{equation} \label{eq:1}
\begin{array}{l}
 i)\ n=3,\ 3<k<5,\\
 ii)\ n=4,\ 2<k<3,\\

 iii)\ n=5,\ \frac{5}{3}<k<\frac{7}{3},\\
 iv)\ n\geq6,\ \frac{n}{n-2}<k\leq\frac{n}{n-3}.\\

\end{array}
\end{equation}
Now we recall a version of the Christ-Kiselev lemma.
\begin{lem}
 Let $X$ and $Y$ be Banach spaces, and for all $s,t\in{\R}^+$ let $K(s,t)$ be an operator from $X$ to $Y$. Suppose that
\[\left\Vert\int_0^{t_0}K(s,t)h(s)ds\right\Vert_{L^l([t_0,+\infty[,Y)}\leq A\Vert h\Vert_{L^r({\R}^+,X)},\]
for some $A>0$, $1\leq r<l\leq+\infty$,   all $t_0\in{\R}^+$ and $h\in L^r({\R}^+,X)$. Then we have

\[\left\Vert\int_0^tK(s,t)h(s) \textrm{d}s\right\Vert_{L^l({\R}^+,Y)}\leq AC_{r,l}\Vert h\Vert_{L^r({\R}^+,X)},\]
where $C_{r,l}>0$ depends only on $r,l$.
 \end{lem}
We refer to  [7] for the proof of Lemma 1 (see also the original paper [1]). Notice that in [7] the above result is formulated with ${\R}$ instead of ${\R}^+$ and $s,t,t_0\in{\R}$, but, as it was mentioned in [7], the same proof works for intervals and in particular for ${\R}^+$. We need the following
\begin{lem}
 Let $a(t,x)$ satisfy the  conditions (1.2). Let $T_1\leq\delta$, and $2\leq p,q<+\infty$ satisfy the conditions (2.1). Then for all $h\in L^1([0,T_1],L^2({\R}^n))$ we have
\begin{equation} \label{eq:1}\left\Vert \int_0^tV(t,s)h(s) \textrm{d}s\right\Vert_{L^p([0,T_1],L^q({\R}^n))}\leq C\Vert h\Vert_{L^1([0,T_1],L^2({\R}^n))}\end{equation}
with $C>0$ independent of $T_1$.
\end{lem}

\begin{proof}
 Let $t_0\in[0,T_1]$. We have
\[\left\Vert\int_0^{t_0}V(t,s)h(s) \textrm{d}s\right\Vert_{L^p([t_0,T_1],L^q({\R}^n))}\leq\int_0^{t_0}\Vert V(t,s)h(s)\Vert_{L^p([t_0,T_1],L^q({\R}^n))} \textrm{d}s.\]
From the definition of $V(t,s)$ we know that
\[\begin{array}{lll}\Vert V(t,s)h(s)\Vert_{L^p([t_0,T_1],L^q({\R}^n))}&=&\left\Vert \left(\mathcal{U}(t,s)(0,h(s))\right)_1\right\Vert_{L^p([t_0,T_1],L^q({\R}^n))}\\
\ &\leq&\left\Vert \left(\mathcal{U}(t,0)(\mathcal{U}(0,s)(0,h(s))\right)_1\right\Vert_{L^p([0,T_1],L^q({\R}^n))}.\end{array}\]
Then, the estimate (2.2) implies that for all $s\in[0,t_0]$ we obtain
\[\begin{array}{lll}\left\Vert \left(\mathcal{U}(t,0)(\mathcal{U}(0,s)(0,h(s))\right)_1\right\Vert_{L^p([0,T_1],L^q({\R}^n))}&\leq& C_\delta\left\Vert\mathcal{U}(0,s)(0,h(s))\right\Vert_{\B}\\
\ &\leq& C'_\delta\Vert h(s)\Vert_{L^2({\R}^n)}, \end{array}\]
where $\displaystyle C'_\delta=C_\delta\sup_{s\in[0,T_1]}\Vert \mathcal U(0,s)\Vert$ is independent of $t_0$.
It follows
\[\begin{array}{lll}\left\Vert\int_0^{t_0}V(t,s)h(s) \textrm{d}s\right\Vert_{L^p([t_0,T_1],L^q({\R}^n))}&\leq&C'_\delta\int_0^{t_0}\Vert h(s)\Vert_{L^2({\R}^n)} \textrm{d}s\\
\ &\leq& C'_\delta\Vert h\Vert_{L^1([0,T_1],L^2({\R}^n))}.\end{array}\]
Consider $K(s,t)=\mathds{1}_{[0,T_1]}(t)\mathds{1}_{[0,T_1]}(s) V(t,s)$, $X=L^2({\R}^n)$ and $Y=L^q({\R}^n)$. Since $p>1$, the Christ-Kiselev lemma yields
\[\left\Vert\int_0^{t}V(t,s)h(s) \textrm{d}s\right\Vert_{L^p([0,T_1],L^q({\R}^n))}\leq C(\delta,p)\Vert h\Vert_{L^1([0,T_1],L^2({\R}^n))}.\]
\end{proof}
 Applying (3.3), we will show  that problem (1.1) is locally well-posed for $k$ and $n$ satisfying  the conditions (3.2).
\begin{Thm}
Assume that $a(t,x)$ is a $C^\infty$ function on ${\R}^{n+1}$ satisfying conditions (1.2) and let $k$ and $n$ satisfy (3.2).  Then there exists  $T_1>0$ such that  problem (1.1) admits a weak solution $u$ on $[0,T_1]$. Moreover, $u$ is the unique weak solution of (1.1) on $[0,T_1]$ satisfying the following properties:
\[\begin{array}{l} (i)\ u\in\mathcal{C}([0,T_1],\dot{H}^1({\R}^n)),\quad 
       (ii)\ u_t\in\mathcal{C}([0,T_1],L^2({\R}^n)),\\
       (iii)\ u\in L^p([0,T_1],L^{2k}({\R}^n)) \quad\textrm{with}\quad \frac{1}{p}=\frac{n(k-1)}{k}-1.\end{array}\]

\end{Thm}
\begin{proof}
Let $k$ and $n$ satisfy (3.2). We have seen  that we can find $2\leq p,q<+\infty$ satisfying conditions  (2.1) so that $\frac{k}{p}<1$ and $\frac{k}{q}=\half$. Consider  the norm $\Vert.\Vert_{Y_{T_1}}$ defined by
\[\Vert u\Vert_{Y_{T1}}=\Vert u\Vert_{\mathcal{C}([0,T_1],\dot{H}^1)}+\Vert u\Vert_{L^p([0,T_1],L^q({\R}^n))}\]
and
\[Y_{T_1}=\mathcal{C}([0,T_1],\dot{H}^1)\bigcap L^p([0,T_1],L^q({\R}^n))\]
with $T_1$ to be determined.
Notice that $(Y_{T_1},\Vert.\Vert_{Y_{T_1}})$ is a Banach space. Assume $f\in{\B}$, $M>0$ and let $B_M=\{u\in Y_{T_1}:\Vert u\Vert_{Y_{T_1}}\leq M\}$, with $M$ to be determined. The problem of finding a weak solution $u$ of (1.1) is equivalent to find a fixed point of the map
\[\mathcal{G}(u)= (\mathcal{U}(t,0)g)_1+\int_0^tV(t,s)f_k(u(s)) \textrm{d}s.\]

Let $u\in B_M$. We have

\[\begin{array}{lll}
 \Vert\int_0^tV(t,s)f_k(u(s)) \textrm{d}s\Vert_{\mathcal{C}([0,T_1],\dot{H}^1)}&\leq&\underset{t\in[0,T_1]}\sup\int_0^{T_1}\mathds{1}_{[0,t]}(s)\Vert  V(t,s)f_k(u(s))\Vert_{\dot{H}^1} \textrm{d}s\\
\ &\leq& \underset{t\in[0,T_1]}\sup\int_0^{T_1}\Vert  V(t,s)f_k(u(s))\Vert_{\dot{H}^1} \textrm{d}s.\\
\end{array}\]
The estimates (2.2)  imply that for $T_1\leq\delta$  there exists $C>0$ independent of $T_1$ such that
\begin{equation}\label{eq:1}\left\Vert\int_0^tV(t,s)f_k(u(s)) \textrm{d}s\right\Vert_{\mathcal{C}([0,T_1],\dot{H}^1)}\leq C\int_0^{T_1}\Vert f_k(u(s))\Vert_{L^2({\R}^n)} \textrm{d}s\leq C_1\int_0^{T_1}\Vert \vert u\vert^k(s)\Vert_{L^2({\R}^n)} \textrm{d}s.\end{equation}
On the other hand, Lemma 2 yields
\begin{equation}\label{eq:1}\left\Vert\int_0^tV(t,s)f_k(u(s)) \textrm{d}s\right\Vert_{L^p([0,T_1],L^q({\R}^n))}\leq C_2\int_0^{T_1}\Vert \vert u\vert^k(s)\Vert_{L^2({\R}^n)} \textrm{d}s.\end{equation}
 We deduce from (3.4) and (3.5)  that
\begin{equation}\label{eq:1}\left\Vert\int_0^tV(t,s)f_k(u(s)) \textrm{d}s\right\Vert_{Y_{T_1}}\leq C_3\int_0^{T_1}\Vert \vert u\vert^k(s)\Vert_{L^2({\R}^n)} \textrm{d}s=C_3\int_0^{T_1}\Vert u(s)\Vert_{L^q}^k \textrm{d}s. \end{equation}
Since $\frac{k}{p}<1$, an application of the   H$\ddot{\textrm{o}}$lder inequality yields
\[\left\Vert\int_0^tV(t,s)f_k(u(s)) \textrm{d}s\right\Vert_{Y_{T_1}}\leq C_3\Vert u\Vert^k_{L^p([0,T_1],L^q({\R}^n))}(T_1)^{1-\frac{k}{p}}\leq C_3M^k(T_1)^{1-\frac{k}{p}}.\]
Let $M$ be such that $\frac{M}{2}\geq 2C(\Vert g_1\Vert_{\dot{H}^1}+\Vert g_2\Vert_{L^2})$ and let $T_1$ be small enough such that
  \[C_3M^k(T_1)^{1-\frac{k}{p}}\leq \frac{M}{2}.\]

  Then  $\Vert \mathcal{G}(u)\Vert_{Y_{T_1}}\leq M$ and $\mathcal{G}(u)\in Y_{T_1}$. We have $\mathcal{G}(B_M)\subset B_M$ and $B_M$ is a closed set of the Banach space $(Y_{T_1},\Vert.\Vert_{Y_{T_1}})$. Now we will  show that we can choose $T_1$ small enough so that  $\mathcal{G}$ becomes a contraction. Let $u,v\in B_M$. We know that
\[\mathcal{G}(u)-\mathcal{G}(v)=\int_0^tV(t,s)(f_k(u(s))-f_k(v(s))) \textrm{d}s.\]
In the same way as in the proof of  inequality (3.4), Theorem 1 and Lemma 2 imply that
\[\Vert \mathcal G(u)-\mathcal G(v)\Vert_{Y_{T_1}}\leq C_4\int_0^{T_1}\Vert f_k(u(s))-f_k(v(s))\Vert_{L^2} \textrm{d}s.\]
On the other hand, $f_k$ satisfies
\[\vert f_k(u)-f_k(v)\vert\leq C_5\vert u-v\vert (\vert u\vert+\vert v\vert)^{k-1}.\]
Consequently,
\[\Vert \mathcal G(u)-\mathcal G(v)\Vert_{Y_{T_1}}\leq  C_6\int_0^{T_1}\Vert \vert u(s)-v(s)\vert(\vert u(s)\vert +\vert v(s)\vert)^{k-1}\Vert_{L^2} \textrm{d}s.\]
Since $\frac{k-1}{q}+\frac{1}{q}=\frac{k}{q}=\half$, by the generalized H$\ddot{\textrm{o}}$lder's inequality, we have
\[\begin{array}{lll}\Vert \vert u(s)-v(s)\vert(\vert u(s)\vert +\vert v(s)\vert)^{k-1}\Vert_{L^2}&\leq& \Vert u-v\Vert_{L^q}\Vert(\vert u(s)\vert +\vert v(s)\vert)^{k-1}\Vert_{L^{\frac{q}{k-1}}}\\
\ &\leq& \Vert u-v\Vert_{L^q}(\Vert u\Vert_{L^q}+\Vert v\Vert_{L^q})^{k-1}.\\ \end{array}\]
This leads to
\[\Vert \mathcal G(u)-\mathcal G(v)\Vert_{Y_{T_1}}\leq C_7\int_0^{T_1}\Vert u(s)-v(s)\Vert_{L^q}(\Vert u(s)\Vert_{L^q}+\Vert v(s)\Vert_{L^q})^{k-1} \textrm{d}s.\]
Applying   H$\ddot{\textrm{o}}$lder's inequality ones more, we find
\[\Vert \mathcal G(u)-\mathcal G(v)\Vert_{Y_{T_1}}\leq C_7(T_1)^{1-\frac{k}{p}}2^{k-1}M^{k-1}\Vert u-v\Vert_{L^p([0,T_1],L^q({\R}^n))}\leq C_7(T_1)^{1-\frac{k}{p}}(2M)^{k-1}\Vert u-v\Vert_{Y_{T_1}}.\]
Thus, if we choose   $T_1$ so that
\[C_7(2M)^{k-1}(T_1)^{1-\frac{k}{p}}<1,\]
  $\mathcal G$ will be a contraction from $B_M$ to $B_M$. Consequently, there exists a unique  $u\in Y_{T_1}$  such that
\[\Vert u\Vert_{Y_{T_1}}\leq M\quad\textrm{and}\quad\mathcal G(u)=u.\]
\end{proof}
\begin{rem}
 In contrast to the case $a=1$ (see [12], [13] and [15]) for our argument we must use homogeneous Strichartz estimates. This restriction leads   to  a solution in the energy space ${\B}$. Moreover, we have more restrictions on the values of $k>1$.
\end{rem}

Since we use  estimates (2.2) to prove Theorem 3, the length $T_1$ of the interval $[0,T_1]$ on which the existence result holds, is majored by the length $\delta$ of the interval on which estimates (2.2) are established. To improve this existence result, in the same way as  in [12], we will apply global estimates in the next section. 

\section{Long time  existence for small initial data}
In this section we assume that $n\geq3$ ,  $a(t,x)$ is $T$-periodic with respect to $t$ and (H1), (H2) are fulfilled. We will use the estimates (2.4) to find solutions of (1.1) defined in $[0,T_1]$, with $T_1$ only depending on $k$, $n$ and $g$. For this purpose we must find $k>1$ such that there exist $2\leq p,q<+\infty$ satisfying (2.3) for which
\begin{equation} \label{eq:2} k=\frac{q}{2},\quad\frac{k}{p}<1.\end{equation}
 Then $k>1$ satisfies (4.1) with $p,q$ satisfying (2.3) if the following conditions are fulfilled
\begin{equation} \label{eq:1} \begin{array}{l}
  i)\ n=3,\ 3<k<5,\\
  ii)\ n=4,\ 2<k<3,\\
 iii)\ n=5,\ \frac{5}{3}<k<\frac{7}{3},\\

 iv)\ n\geq6,\ \frac{n}{n-2}<k<\frac{n}{n-3}.\end{array}\end{equation}

\begin{lem}
Assume that (H1) and (H2) are fulfilled, $a(t,x)$ is $T$-periodic with respect to $t$ and $n\geq3$.  Let $t\geq s\geq0$. Then
\[\Vert \mathcal{U}(t,s)\Vert_{\mathcal{L}({\B})}\leq C_0\]
with $C_0>0$ independent of $s$ and $t$.
\end{lem}
\begin{proof}
 Let $m\in\mathbb{N}$ be such that $0\leq s-mT<T$. We have
\[\mathcal{U}(t,s)=\mathcal{U}(t-mT,s-mT)=\mathcal{U}(t-mT,0)\mathcal{U}(0,s-mT).\]
Since $t-mT\geq s-mT\geq0$, Theorem 2 implies
\[\Vert \mathcal{U}(t-mT,0)\Vert_{\mathcal{L}({\B})}\leq C'\]
with $C'>0$ independent of $t$. Also we have
\[\Vert \mathcal{U}(0,s-mT)\Vert_{\mathcal{L}({\B})}\leq \sup_{s'\in[0,T]}\Vert \mathcal{U}(0,s')\Vert_{\mathcal{L}({\B})}= C''.\]
It follows  that
\[\Vert \mathcal{U}(t,s)\Vert_{\mathcal{L}({\B})}\leq C'C''=C_0\]
and $C_0$ is independent of $t$ and $s$.

\end{proof}

The estimates (2.4), the Christ-Kiselev lemma and Lemma 3 imply the following
\begin{lem}
 Assume  $n\geq3$  and let $a(t,x)$  be  $T$-periodic with respect to $t$ such that (H1) and (H2) are fulfilled. Let $2\leq p,q<+\infty$ satisfy condition (2.3) and let $T_1>0$. Then for all \\
$h\in L^1([0,T_1],L^2({\R}^n))$ we have
\[\left\Vert\int_0^tV(t,s)h(s) \textrm{d}s\right\Vert_{L^p([0,T_1],L^q({\R}^n))}\leq C\Vert h\Vert_{L^1([0,T_1],L^2({\R}^n))}\]
with $C>0$ independent of $g$ and $T_1$.
\end{lem}

\begin{proof}
Let $t_0>0$ , $s\in[0,t_0]$ and $t>t_0$. Consider $mT\leq t_0<(m+1)T$. We have
\[\begin{array}{lll}V(t,s)h(s)&=&\left(\mathcal{U}(t-mT,s-mT)( 0,h(s))\right)_1\\
                                                                  
\ &=&\left(\mathcal{U}(t-mT,0)\mathcal{U}(0,s-mT)(0,h(s))\right)_1\\
 \ &=&\left(\mathcal{U}(t-mT,0)\mathcal{U}(mT,s)(0,h(s))\right)_1.
                                                                   \end{array}\]
Thus, the estimate (2.4) implies
\[\left\Vert \int_0^{t_0}V(t,s)h(s) \textrm{d}s\right\Vert_{L^p([t_0,+\infty[,L^q({\R}^n))}\leq C\int_0^{t_0}\left\Vert \mathcal{U}(mT,s)(0,h(s))\right\Vert_{\B} \textrm{d}s.\]
Since  $0\leq s\leq mT$ for $s\in [0,mT]$, Lemma 3 yields
\[\int_0^{mT}\left\Vert \mathcal{U}(mT,s)(0,h(s))\right\Vert_{\B} \textrm{d}s\leq C_0\int_0^{mT}\Vert h(s)\Vert_{L^2({\R}^n)} \textrm{d}s\]
with $C>0$ independent of $t_0$. In the same way, since $mT\leq t_0<(m+1)T$ we have
\[\int_{mT}^{t_0}\left\Vert \mathcal{U}(0,s-mT)(0,h(s))\right\Vert_{\B} \textrm{d}s\leq \sup_{s\in[0,T]}\Vert \mathcal{U}(0,s)\Vert \int_{mT}^{t_0}\Vert h(s)\Vert_{L^2({\R}^n)} \textrm{d}s.\]
It follows that
\[\left\Vert \int_0^{t_0}V(t,s)h(s) \textrm{d}s\right\Vert_{L^p([t_0,+\infty[,L^q({\R}^n))}\leq C\int_0^{+\infty}\Vert h(s)\Vert_{L^2({\R}^n)} \textrm{d}s.\]
Since $p>1$, the Christ- Kiselev lemma implies
\[\left\Vert \int_0^t V(t,s)h(s) \textrm{d}s\right\Vert_{L^p({\R}^+,L^q({\R}^n))}\leq C_p\Vert h\Vert_{L^1({\R}^+,L^2({\R}^n))}.\]
We deduce that
\[\begin{array}{lll}\left\Vert \int_0^tV(t,s)h(s) \textrm{d}s\right\Vert_{L^p([0,T_1],L^q({\R}^n))}&\leq&\left\Vert \int_0^tV(t,s)\mathds{1}_{[0,T_1]}(s)h(s) \textrm{d}s\right\Vert_{L^p({\R}^+,L^q({\R}^n))}\\

\ &\leq& C\Vert h\Vert_{L^1([0,T_1],L^2({\R}^n))}\end{array}\]
with $C>0$  independent of $T_1$.
\end{proof}

\begin{Thm}
Assume that $k$ and $n$ satisfy  the  conditions (4.2). Let $a(t,x)$ be  $T$-periodic with respect to $t$ and let (H1), (H2) be fulfilled. Then  there exists $C(k,f_k,T,\rho,n)$ such that for all $g\in{\B}$ we can find a weak solution $u$ of (1.1) on $[0,T_1]$ with
\begin{equation} \label{eq:1}T_1=C(k,f_k)\left(\Vert g\Vert_{\B}\right)^{-d},\end{equation}
where $\displaystyle d=\frac{2(k-1)}{(n+2)-(n-2)k}$. Moreover, $u$ is the unique weak solution of (1.1) on $[0,T_1]$ satisfying the following properties:
\begin{equation} \label{eq:1}\begin{array}{l} (i)\ u\in\mathcal{C}([0,T_1],\dot{H}^1({\R}^n)),\quad
       (ii)\ u_t\in\mathcal{C}([0,T_1],L^2({\R}^n)),\\
       (iii)\ u\in L^p([0,T_1],L^{2k}({\R}^n)) \quad\textrm{with}\quad \frac{1}{p}=\frac{n(k-1)}{k}-1.\end{array}\end{equation}

\end{Thm}
\begin{proof}
Let $C_f>0$ be such that $\vert f_k(u)\vert\leq C_f\vert u\vert^k$ and $\vert f_k(u)-f_k(v)\vert\leq C_f \vert u-v\vert\left(\vert u\vert+\vert v\vert\right)^{k-1}$. Then, Theorem 2 and Lemma 4 imply that there exists $A_k$ such that for all $T_1>0$
\[\left\Vert\int_0^tV(t,s)h(s) \textrm{d}s\right\Vert_{L^p([0,T_1],L^q({\R}^n))}\leq A_k\Vert h\Vert_{L^1([0,T_1],L^2({\R}^n))}\]
and
\[\left\Vert(\mathcal{U}(t,0)g)_1\right\Vert_{L^p([0,T_1],L^q({\R}^n))}\leq A_k\Vert g\Vert_{\B}.\]
According to the proof of Theorem 3, $\mathcal{G}(u)= (\mathcal{U}(t,0)g)_1+\int_0^tV(t,s)f_k(u(s)) \textrm{d}s$
admits a fixed point  in the set \[\{u\in\mathcal{C}([0,T_1],\dot{H}^1)\cap L^p([0,T_1],L^q):\ \Vert u\Vert_{\mathcal{C}([0,T_1],\dot{H}^1)}+\Vert u\Vert_{ L^p([0,T_1],L^q)}\leq M\}\] if we choose $M, T_1>0$ so that
\begin{equation} \label{eq:1}\left\{ \begin{array}{l}\displaystyle A_k\Vert g\Vert_{\B}+C_3M^k(T_1)^{1-\frac{k}{p}}\leq M,\\
\displaystyle  C_7(2M)^{k-1}(T_1)^{1-\frac{k}{p}}<1.\\
\end{array}\right.
\end{equation}
In particular (4.5) will be fulfilled if
\begin{equation} \label{eq:1}\left\{ \begin{array}{l}\displaystyle A_k\Vert g\Vert_{\B}+C_3M^k(T_1)^{1-\frac{k}{p}}= M,\\
\displaystyle  C_7(2M)^{k-1}(T_1)^{1-\frac{k}{p}}<1.\\
\end{array}\right.
\end{equation}
We will choose $M,T_1$ so that (4.6) holds. Let $t_1=(T_1)^{1-\frac{k}{p}}$. We find that the system (4.6) is equivalent to the following

\begin{equation}\label{eq:1}\left\{ \begin{array}{l}\displaystyle t_1=\frac{M-A_k\Vert g\Vert_{\B}}{C_3M^k},\\
\ \\
 \displaystyle0<\frac{M-A_k\Vert g\Vert_{\B}}{M}<\frac{1}{C_72^{k-1}}.\\
\end{array}\right.\end{equation}
Since $M\longmapsto\frac{M-A_k\Vert g\Vert_{\B}}{M}$ is strictly increasing, we obtain that $(t_1,M)$ is a solution of (4.7) if
\[M<\frac{2^{k-1}A_k\Vert g\Vert_{\B}}{2^{k-1}-1}.\]
Take 
\[M_0=\frac{\alpha2^{k-1}A_k\Vert g\Vert_{\B}}{2^{k-1}-1}\quad\textrm{and}\quad t_1=\frac{M_0-A_k\Vert g\Vert_{\B}}{C_3(M_0)^k}\]
with $\displaystyle 1-\frac{1}{2^{k-1}}<\alpha<1$.
Then $(M_0,t_1)$ is a solution of (4.7) and we have

\[t_1=\frac{\frac{\alpha2^{k-1}}{2^{k-1}-1}-1}{C_8\left(\frac{2^{k-2}A_k}{2^{k-1}-1}\right)^k\Vert g\Vert_{\B}^{k-1}}=C'(k,f_k)\Vert g\Vert_{\B}^{-(k-1)}.\]

Thus for $M=M_0$ and  
\[T_1=(t_1)^{\frac{1}{1-\frac{k}{p}}}=C(k,f_k)\left(\Vert g\Vert_{\B}\right)^{-\frac{k-1}{1-\frac{k}{p}}},\]
$M$ and $T_1$ satisfy conditions (4.6).
Moreover, we know that
\[\frac{k}{p}=\frac{(n-2)k-n}{2}.\]
Thus we have
\[\frac{k-1}{1-\frac{k}{p}}=\frac{2(k-1)}{(n+2)-(n-2)k}\]
and $M$, $T_1$ satisfy conditions (4.5) if $M=M_0$ and
\[T_1=C(k,f_k)\left(\Vert g\Vert_{\B}\right)^{-\frac{2(k-1)}{(n+2)-(n-2)k}}.\]  
Note that for $n\geq6$ we have
$\displaystyle\frac{n}{n-3}\leq\frac{n+2}{n-2}$ and $\displaystyle k<\frac{n}{n-3}$ leads to $\displaystyle k<\frac{n+2}{n-2}$.

\end{proof}

\begin{rem}
Let $\Vert g\Vert_{\B}=\epsilon$  and $T_1=C\epsilon^{-d}$, $C,d>0$ being the constants defined by (4.3). Then Theorem 4 implies that there exists a unique   solution of (1.1) satisfying (4.4).  
\end{rem}

\end{document}